\documentclass[11pt]{article}
\usepackage{amsfonts,amsmath,amssymb}
\usepackage{latexsym}
\usepackage{amsthm}
\usepackage{amscd}
\usepackage{commath}
\usepackage{epsfig}
\usepackage{amssymb}
\usepackage{caption}
\usepackage{enumitem}
\usepackage{mathtools}
\usepackage{booktabs}
\usepackage[table]{xcolor}
\usepackage[toc,page]{appendix}
\usepackage{array}
\usepackage{stmaryrd}
\newcolumntype{P}[1]{>{\centering\arraybackslash}p{#1}}
\usepackage{siunitx}

\usepackage[english]{babel}

\usepackage[letterpaper,top=2cm,bottom=2cm,left=2.5cm,right=2.5cm,marginparwidth=1.5cm]{geometry}

\usepackage[colorlinks=true, allcolors=blue,bookmarksnumbered]{hyperref}
\hypersetup{
	citecolor = red!80!black
}
\usepackage{color,soul}
\setstcolor{red}
\usepackage{float,subcaption}
\usepackage{graphicx}
\usepackage{stackengine}
\usepackage{tikz}
\usetikzlibrary{intersections}
\usepackage[export]{adjustbox}
\usepackage{array}
\RequirePackage{aliascnt}

\usepackage[capitalise]{cleveref}

\usepackage{microtype}

\usepackage[percent]{overpic}

\usepackage{lineno}

\newtheorem{defn}{Definition}[section]

\newaliascnt{proposition}{defn}
\newtheorem{prop}[proposition]{Proposition}
\aliascntresetthe{proposition}
\newaliascnt{lemma}{defn}
\newtheorem{lem}[lemma]{Lemma}
\aliascntresetthe{lemma}
\newaliascnt{theorem}{defn}

\newtheorem{theorem}[theorem]{Theorem}
\aliascntresetthe{theorem}

\newtheorem{cor}[defn]{Corollary}
\theoremstyle{remark}
\newtheorem{remark}[defn]{Remark}

\numberwithin{equation}{section}
\numberwithin{figure}{section}

\newcommand{\bb}{\begin{equation}}
\newcommand{\ee}{\end{equation}}

\newcolumntype{?}{!{\vrule width 1.5pt }}
\newlength\savedwidth

\newcommand{\tightoverset}[2]{%
  \mathop{#2}\limits^{\vbox to -.5ex{\kern-0.75ex\hbox{$#1$}\vss}}}

\newcommand{\dfn}[1]{\textbf{\textit{#1}}}
\newcommand\Leb{\lambda}
\newcommand{\R}{\mathbb{R}}

\def\rproof{\begin{proof}}

\def\Qed{\end{proof}}

\usepackage{bbm}
\usepackage{xspace}
\usetikzlibrary{arrows.meta}

\newcommand\numberthis{\addtocounter{equation}{1}\tag{\theequation}}

\newcommand{\wrt}{\mbox{w.r.t.}\xspace}
\newcommand{\iid}{\mbox{i.i.d.}\xspace}

\newcommand{\ie}{\mbox{i.e.}\xspace}
\newcommand{\eg}{\mbox{e.g.}\xspace}

\newcommand{\inv}{^{-1}}
\newcommand{\invv}{^{-2}}

\newcommand{\oneton}[1]{#1_1,\dots,#1_n}

\newcommand{\fallfac}[2]{\prescript{}{#1}{\!\left(#2\right)}}

\newcommand{\ind}{\mathbbm1}
\newcommand{\bfo}[1][]
{\ind
   \ifthenelse{\equal {#1} {}}
    {}{
        _{\{#1\}}
    }
}
\DeclareMathOperator{\MAD}{MAD}

\newcommand{\nocontentsline}[3]{}
\newcommand{\tocless}[2]{\bgroup\let\addcontentsline=\nocontentsline#1{#2}\egroup}

\def\N{\mathbb N}

\def\R{\mathbb R}

\newcommand{\Nn}{\mathbb N_{\geq0}}

\newcommand{\PR}{\mathbb P}
\renewcommand{\P}{\mathbb P}

\newcommand{\OO}{\mathcal O}

\NewDocumentCommand\E{o o}
{\mathbb E
   \IfNoValueTF{#2}{}{_{#2}}
    \IfNoValueTF {#1}
   {}
    {\!\left[#1\right]\!}
}

\newcommand{\sigA}{\mathcal A}

\newcommand{\oo}{\mathbf o}

\DeclareMathOperator{\Var}{Var}
\DeclareMathOperator{\Cov}{Cov}

\newcommand{\LebR}{\Leb_\R}

\newcommand{\dd}{\mathop{}\!\mathrm{d}}
\newcommand{\intnsimplex}[1][]
{\ifthenelse{\equal {#1} {}}
    {
    \int_{\R^{n+1}_+}\bfo_{\{\sum_{i=1}^{n+1}x_i=1\}}
    }
    {
    \int_{\R^{n+1}_+}\bfo_{\{\sum_{i=1}^{n+1}x_i=#1\}}
    }
}

\newcommand{\ltwonorm}[1]{\lVert#1\rVert_2}

\newcommand{\Sn}{\mathcal S_n}

\newcommand{\QM}[4]{{\rm QM}_{#4}(\mathbf{#1},\mathbf{#2},#3)}

\newcommand{\oinO}{\omega\in\Omega}

\newcommand{\cdf}{cumulative distribution function\xspace}

\def\probspaceA{$(\Omega,\sigA,\PR)$\xspace}

\newcommand{\thispropname}{No long cycles}

\newcommand{\CNES}{Centre national d'\'etudes spatiales\xspace}

\title{\centering\textsc{The Poissonian limit of a model inspired by the one-dimensional \scalebox{0.8}{ERAP}, and its application to satellite data}}
\author{
Matteo \textsc{D'Achille}\thanks{Institut Élie Cartan de Lorraine, CNRS, Universit\'e de Lorraine,  F-57070, Metz, France.\newline  $_{}$ \hfill \href{mailto:matteo.d-achille@univ-lorraine.fr}{\texttt{matteo.d-achille@univ-lorraine.fr}}}, Jessie \textsc{Levillain}\thanks{Centre National d’\'Etudes Spatiales, 31401 Toulouse; INRAE MaIAGE, Université Paris-Saclay, 78352 Jouy-en-Josas CEDEX, France.\hfill  \href{mailto:jessie.levillain@inrae.fr}{\texttt{jessie.levillain@inrae.fr}}}, Pascal \textsc{Moyal}\thanks{Laboratoire de mathématiques de Versailles, CNRS, Université de Versailles Saint-Quentin-en-Yvelines, F-78000 Versailles, France. Previously at: Institut Élie Cartan de Lorraine, CNRS, Universit\'e de Lorraine, \& LORIA, F-54000, Nancy, France.\hfill  \href{mailto:pascal.moyal@uvsq.fr}{\texttt{pascal.moyal@uvsq.fr}}}, Achim \textsc{Schlather}\thanks{Institut Élie Cartan de Lorraine, CNRS, Universit\'e de Lorraine,  F-57070, Metz; Inria, LORIA, F-54000, Nancy, France. \hfill  \href{mailto:achim.schlather@univ-lorraine.fr}{\texttt{achim.schlather@univ-lorraine.fr}}}
}
\date{\today}

\begin{document}
\maketitle
\begin{abstract}
We consider a model of particle tracking in one dimension studied by Quarck in which $n$ independent random points of the interval $[0,1]$ are shifted by \iid noise at a common scale $\varepsilon_n>0$. We prove that, if $\varepsilon_n=\alpha n\invv+\oo(n\invv)$ for some $\alpha>0$ when $n\to \infty$, then in the latter limit the total number of consecutive swaps in the shifted sample \wrt its order statistics converges in distribution to a Poisson random variable with a parameter depending explicitly on the product of $\alpha$, the mean absolute difference of the shifts and the $L^2$-norm of the density of points. We compare the statistical criterion implied by this theorem to numerical simulations based on real terrain maps, and discuss both its usefulness and limitations.
\end{abstract}

\section{Introduction and model definition}
\paragraph{Satellite data acquisition.}
When a satellite flies over the surface of some planet, one might be interested in identifying certain objects, \eg bushes, trees, craters or edges of bodies of water \cite{gasnier}, in consecutive images and keeping track of which object is which. However, since the satellite moves between images, the objects in the second image are shifted compared to the first image. To model this scenario, we assume that all $n$ objects of interest are indistinguishable and are moved independently and randomly between images. This is different from other approaches with  fixed objects that allow the utilization of other data \cite{gasnierDatabase,geo} and also disallows classification-based approaches. Further, we assume that no new objects appear and that no objects vanish.

\paragraph{Particle tracking and the Euclidean Random Assignment Problem.} 
The problem of particle tracking is classical and has been addressed from many angles, especially in statistical physics \cite{CKZ10}.

On the random geometry side, a natural toy model is provided by the \dfn{Euclidean Random Assignment Problem (ERAP)}, see~\cite{dachille}. This is defined as the study of the random optimal energy $\mathcal{H}_{\rm opt}$ of assigning $n$ points (blues) to $n$ points (reds) randomly distributed on a Polish metric space $(\Omega,D)$ depending on the choice of $(\Omega,D)$, on the law(s) of the points and on an energy-distance exponent $p \in \mathbb{R}$. When $\Omega=[0,1]^d$ is endowed with $D=$ the Euclidean distance, it is known that, for $d\geq 3$, $\E[\mathcal{H}_{\rm opt}]\sim c_{p,d} n^{1-p/d}$ as $n \to \infty$ for some constant $c_{p,d}$ independent of $n$; and that $\mathcal{H}_{\rm opt}$ concentrates in this limit as long as $p\in [1,d^2/2)$, see \cite[Theorem 2.1]{DMT26}. Based on a $p\to q$ transfer principle and extensive numerical simulations, it has been conjectured in \cite[Conjecture 1.1]{DMT26} that concentration holds for all $p\geq 1$ in any dimension $d\geq 3$. This conjecture has been recently settled in \cite{GKT26}. When $\Omega=[0,1]^2$ (resp.~a closed 2-dimensional Riemannian manifold) is endowed with the Euclidean distance (resp.~the metric of the manifold), the problem is much less understood. Addressing certain aspects of a celebrated linearization prediction by Caracciolo--Lucibello--Parisi--Sicuro~\cite{CLPS}, Ambrosio--Stra--Trevisan used PDE methods to prove in~\cite[Theorem 1.1]{AST19} that, for the quadratic case $p=2$, as $n\to \infty$  
\begin{equation}\label{eq.2dexp}
\E[\mathcal{H}_{\rm opt}] = \frac{{\rm Vol}(\Omega)}{2\pi} \ln(n) + \oo(\ln n)\; ,
\end{equation}
where ${\rm Vol}(\Omega)$ denotes the volume of $\Omega$. The $\oo(\ln n)$ appearing in \Cref{eq.2dexp} was later improved to $\oo(\sqrt{\ln(n) \ln \ln (n)})$ by Ambrosio--Glaudo in~\cite[Theorem 1.2]{AG19}. It is conjectured in~\cite{BCCDSS21} that the constant part of the remainder term in \cref{eq.2dexp} is given by the Robin mass of the Laplace--Beltrami operator on  $\Omega$ (or equivalently by its Kronecker mass, due to a celebrated result of Morpurgo~\cite{Mo02}). To the best of the authors' knowledge, whether $\mathcal{H}_{\rm opt}$ concentrates in $d=2$ (say, on the unit square or torus) is an open problem.
 
\begin{figure}[htb]
\centering
    \begin{tikzpicture}[
    >={Stealth},
    every node/.style={circle, minimum size=10mm},x=1.2cm
]

\node[fill=blue,text=white] (A) at (-4,0) {$b_{(1)}$};
\node[fill=blue,text=white] (B) at (3,0) {$b_{(2)}$};
\node[fill=blue,text=white] (C) at (4,0) {$b_{(3)}$};
\node[fill=blue,text=white] (D) at (6,0) {$b_{(4)}$};
\node[fill=red] (E) at (-2,0) {$r_1$};
\node[fill=red] (F) at (7,0) {$r_3$};
\node[fill=red] (G) at (1,0) {$r_2$};
\node[fill=red] (H) at (2,0) {$r_4$};
\draw[->] (A) to[bend right=30] (E);
\draw[->] (B) to[bend right=50] (G);
\draw[->] (C) to[bend right=40] (F);
\draw[->] (D) to[bend right=40] (H);
\draw[|->,draw=purple,shorten <=1pt] (E) to[bend right=30] (A);
\draw[|->,draw=purple,shorten <=1pt] (G) to[bend right=50] (B);
\draw[|->,draw=purple,shorten <=1pt] (H) to[bend right=50] (C);
\draw[|->,draw=purple,shorten <=1pt] (F) to[bend right=60] (D);
\end{tikzpicture}
\caption{Realization of a Quarck model with $n=4$. The arrows $\rightarrow$ signify the addition of some random shift $Z_1,\dots,Z_4$ to each of the blue points. The arrows ${\color{purple}\mapsto}$ represent the solution of the ERAP.}\label{fig.def.model}
\end{figure}
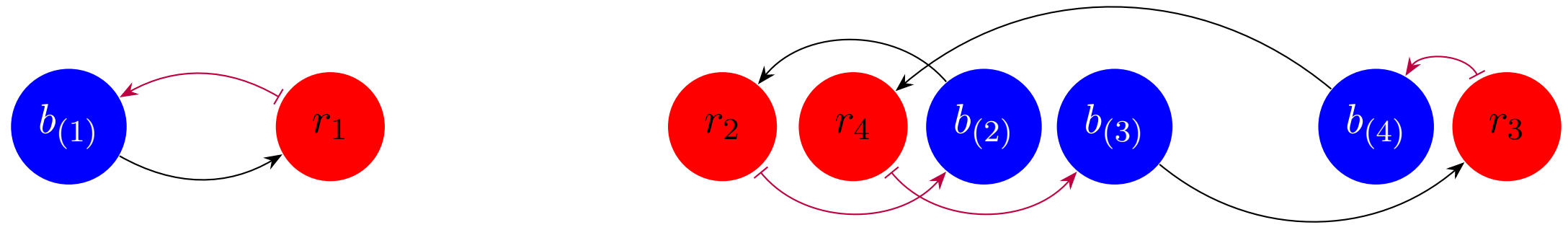

\paragraph{Model.} We simplify the situation to the one-dimensional case by projecting the objects onto a line orthogonally to the flight path. If, with this simplified approach, an identification between images is possible, it also is for the original 2-dimensional case.
Borrowing from previous work of Quarck \cite{quarck}, we formalize this idea as a \dfn{Quarck model}. For the identification in real pictures, one does not have full information about the shift distribution, and possibly neither about the distribution of the blue points. In this model, the goal is, given a set of blue points and a same-sized set of red points, to find the association $B_{(i)}\mapsto R_i$. $B_{(i)}$ denotes the usual order statistic, \ie $B_{(i)}$ is the value of the $i$-th smallest blue point among $n$. One way of trying to obtain this is to view a realization of a Quarck model as a one-dimensional ERAP. In this one-dimensional case, for a convex and continuous cost function $c$, one can show using the Hardy-Littlewood-P\'olya rearrangement inequality \cite{HLP} that
\[\sum_{i=1}^nc(b_{(i)}-r_{(i)})\leq\sum_{i=1}^nc(b_{(i)}-r_{(\sigma(i))})\]
for any permutation $\sigma$ of the symmetric group $\Sn$ over $\llbracket1,n\rrbracket$.  In other words, the solution of the ERAP for a Quarck model is given by $R_{(i)}\mapsto B_{(i)}$ (\cref{fig.def.model}).
In practice, the general ERAP is known to be $P$ complete and can be solved with the \emph{Hungarian method} \cite{kuhn1955,Munkres1957}, which was later improved to be of cubic time complexity \cite{edmonds1972}. In our one-dimensional case, the time complexity to solve the ERAP is $\OO(n\log n)$, as only sorting is required. We use the Bachmann-Landau notation: Let $f,g$ be functions $:A\subseteq\R\rightarrow\R$. Then, we write $f(x)=\OO(g(x))$ if there exists $\alpha>0$ such that for all $x\in A:|f(x)|<\alpha g(x)$, and 
$f(x)=\oo(g(x))$ if $\lim_{x\rightarrow\infty}f(x)/g(x)=0$.
\paragraph{}
We now properly define the model at hand. Throughout, we fix a probability space \probspaceA.
\begin{defn}[Quarck model]
    Let \[\oneton B:\Omega\rightarrow[0,1]\] be \iid random variables \wrt a probability density $f_B$. Further, let $Z_1,\dots,Z_n$ be \iid random variables {with a common \cdf $F_Z$}, also independent from $\oneton B$. We define, for $\varepsilon>0$,
    \[R_i\coloneq B_{(i)}+\varepsilon Z_i.\]
    We call the $B_{(i)}$'s `blue' points, the $R_i$'s `red' points and the $Z_i$ `shifts'. An example is shown in \cref{fig.def.model}.
    For $n\in \N$, we shall denote the Quarck model corresponding to blue points $\mathbf{B}=(B_1,\ldots,B_n)$ and shifts $\mathbf{Z}=(Z_1,\ldots,Z_n)$ by $\QM{B}{Z}{\varepsilon}{n}$.
\end{defn}

\section{Poissonian limit for the error count}
\label{sec:main}
Differently from the usual ERAP, we are not interested in the cost of the assignment but instead in the number of errors from the `correct' assignment, 
$B_{(i)}\mapsto R_i,i\in\llbracket 1,n \rrbracket$. 
Since the solution of the ERAP is sensitive to ordering, it motivates quantifying errors by the number of instances $i\in\llbracket 1,n-1 \rrbracket$ where $R_i>R_{i+1}$. 
\begin{defn}[Error count]\label{def.errorCount}
    We call the random variable
    \[V_n\coloneq \sum_{i=1}^{n-1}\bfo[R_i>R_{i+1}]\]
    the \emph{error count} of the  Quarck model $\QM{B}{Z}{\varepsilon}{n}$.
\end{defn}
In this work, we focus on the distributional limit of $V_n$ when $n\rightarrow\infty$ under some natural assumptions on the law of the blue points and of the shifts.
 
As the number of points is typically large, we will let $n$ increase to infinity, and as is natural we will scale the order of the shifts with $n$, and denote it 
as $\varepsilon_n$. 
As will be shown below, a phase transition occurs for the number of exchanges $V_n$ as a function of $\varepsilon_n$. What we call the {\em critical scale} is the order 
of magnitude for $\varepsilon_n$ at which this transition occurs. It is properly defined below.
\begin{defn}
We say that the Quarck model $\QM{B}{Z}{\varepsilon_n}{n}$ is {\em at the critical scale} if $$\varepsilon_n=\alpha n\invv+\oo(n\invv),\quad n\ge 1.$$
\end{defn}

\noindent Our main result is the following:
\begin{theorem}[Explicit Poissonian limit of the error count]
\label{thm:main}
Consider a Quarck model $\QM{B}{Z}{\varepsilon_n}{n}$ such that 
\begin{itemize}
\item $f_B$ has finite $L^2$-norm w.r.t. Lebesgue measure; 
\item The shifts $Z_1,\dots,Z_n$ have finite $k$-th moment, for some $k>3/2$.
\end{itemize}
Then, at the critical scale we have the convergence in distribution 
$$V_n \overset{\mathcal{D}}{\underset{n\to \infty}{\longrightarrow}} {\rm Poiss}(\mathfrak{c}_\alpha) \, ,$$
where $$\mathfrak{c}_\alpha\coloneqq \frac{\alpha}{2}\, \E[|Z_1-Z_2|]\lVert f_B\rVert_2^2.$$
\end{theorem}
\noindent The proof of Theorem \ref{thm:main} is deferred to Section \ref{sec:proofthm}. 

We use the Copernicus GLO-30 Digital Elevation Model \cite{data} to compare the theoretical model with data artificially generated based on real terrain maps. We find that the Quarck model seems to provide a worst-case scenario as real terrain is locally of similar height and not of independent random heights (see Table \ref{tab:highlowredrawn}).

\section{Preliminary results}
\label{sec:prelimres}
We denote by $\LebR$ Lebesgue measure on $\R$. 
For a probability density $f$ on $\R$, the $L^p$ norm of $f$ is understood w.r.t. Lebesgue measure, that is,
    \[\norm f_p\coloneq\bigg(\int_\R f^p(x)\dd\LebR(x)\bigg)^{1/p},\quad\text{for }p\geq1.\]
Throughout the section, fix a Quarck model $\QM{B}{Z}{\varepsilon_n}{n}$ such that 
$f_B$ has finite $L^2$-norm and 
the shifts have finite $k$-th moment, for some $k>3/2$, and are not almost surely constant (otherwise $V_n$ is trivially 0 for all $n$). 

\subsection{On the \texorpdfstring{$L^2$}{L2} norm of the density of blue points}
{First observe the following convenient characterization of the \texorpdfstring{$L^2$}{L2} norm of the density of blue points as a limit.} 
\begin{lem}\label{lemma.intfxfx+y->ltwonorm_f}
    Let $f$ be a probability density on $[0,1]$ with $\ltwonorm f<\infty$. Then we have the convergence
        \[\int_{[0,{1-y}]}f(x)f(x+y)\dd\LebR(x)\xrightarrow{y\rightarrow\infty}\ltwonorm f^2,\]
    from below. 
\end{lem}
\begin{proof}
    First, observe that, for $y>0$,
    \begin{align*}
        2\int_{[0,{1-y}]}\!\!f(x)f(x+y)\dd\LebR(x)&=\int_{[0,{1-y}]}\!\!f^2(x+y)\dd\LebR(x)+\int_{[0,{1-y}]}\!\!f^2(x)\dd\LebR(x)\\
        &\qquad-\int_{[0,{1-y}]}\!(f(x+y)-f(x))^2\dd\LebR(x).
    \end{align*}
    Let for all $x\in\R$, $\tilde f(x)=f(x)\mathbf 1_{[0,1]}(x).$ Then, $\ltwonorm{\tilde f(\cdot+y)-\tilde f}$ converges to zero as $y\rightarrow0$ \cite[Proposition 8.5]{folland1999real}. That gives
    \begin{align*}
        0\leq\int_{[0,{1-y}]}(f(x+y)-f(x))^2\dd\LebR(x)&\leq\int_{[0,1]}(\tilde f(x+y)-\tilde f(x))^2\dd\LebR(x)\\&=\ltwonorm{\tilde f(\cdot+y)-\tilde f}^2\xrightarrow{y\rightarrow0} 0.
    \end{align*}
    Further,
    \begin{align*}
        \int_{[0,{1-y}]}\!\!f^2(x+y)\dd\LebR(x)+\int_{[0,{1-y}]}\!\!f^2(x)\dd\LebR(x)=&\int_{[y,1]}\!f^2(x)\dd\LebR(x)+\int_{[0,{1-y}]}\!\!f^2(x)\dd\LebR(x)\\
        \xrightarrow{y\rightarrow0}&\,2\ltwonorm f^2,
    \end{align*}
    implying that
   \[\int_{[0,{1-y}]}\!\!f(x)f(x+y)\dd\LebR(x)\xrightarrow{y\rightarrow\infty}\ltwonorm f^2\]
   converges from below.
\end{proof}

\subsection{On the mean absolute differences of the shifts}
{We now introduce a few technical results to characterize $\E\left[|Z_1-Z_2|\right]$, which will be useful in the proofs. Recall that $F_Z$ denotes the common \cdf of 
the shifts. 
\begin{defn}
    We define for $\Delta\geq0$, 
    \[\varphi_Z(\Delta)\coloneq\int_\R F_Z(y-\Delta)\dd\mu_Z(y),\]
     where $\mu_Z((a,b])\coloneq F_Z(b)-F_Z(a)$ is uniquely extended to a measure using Carathéodory's extension theorem \cite[Theorem 1.14]{folland1999real}.
\end{defn}
\noindent Importantly enough, notice that 
\begin{equation}
\label{eq:interpretvarphiZ}
\varphi_Z(\Delta)=\P(Z_1-Z_2 \ge \Delta)=\P(Z_1-Z_2 \le -\Delta), \mbox{ for all } \Delta\ge 0.
\end{equation}
}

\begin{lem}\label{lemma.MAD}    
    {Under the ongoing assumptions we have that 
    \[\int_{\R_{\geq0}}\varphi_Z(\Delta)\dd\LebR(\Delta)=\frac12\E[|Z_1-Z_2|].\]}
\end{lem}
\begin{proof}
    First we use Tonelli's theorem:
    \[\int_{\R_{\geq0}}\varphi_Z(\Delta)\dd\LebR(\Delta)=\int_{\R_{\geq0}}\int_\R F_Z(y-\Delta)\dd\mu_Z(y)\dd\LebR(\Delta)=\int_\R \int_{\R_{\geq0}}\!F_Z(y-\Delta)\dd\LebR(\Delta)\dd\mu_Z(y).\]
    We calculate, using integration by parts for absolutely continuous functions together with the fact that $\lim_{x\rightarrow-\infty}xF_Z(x)=0$,
    \begin{align*}
        \int_\R \int_{\R_{\geq0}}\! F_Z(y-\Delta)\dd\LebR(\Delta)\dd\mu_Z(y)
        &=\int_\R \left.\Delta F_Z(y-\Delta)\right|_{\Delta=0}^{\Delta\rightarrow\infty}\dd\mu_Z(y)+\int_\R\int_{\R_{\geq0}}\! \Delta\dd\mu_{Z}(y-\Delta)\dd\mu_Z(y)\\
        &=\int_\R\int_{\R_{\geq0}}\! \Delta\dd\mu_{Z}(y-\Delta)\dd\mu_Z(y)\\
        &=\int_\R\int_{(-\infty,y]} (y-u)\dd\mu_Z(u)\dd\mu_Z(y)\\
        &=\frac12\int_\R\int_\R|y-u|\dd\mu_Z(u)\dd\mu_Z(y)\\
        &=\frac12\E[|Z_1-Z_2|].
    \end{align*} 
\end{proof}

\subsection{Simplicity of the Quarck model}

\begin{defn}[Assignment in a Quarck model]
The random permutation $\Pi\in\Sn$ given by
    \[\Pi(i)=j\Longleftrightarrow R_i=R_{(j)}\]
    is called the \emph{assignment} associated to a Quarck model $\QM{B}{Z}{\varepsilon}{n}$. 
\end{defn}
{Observe that, as a consequence of \eqref{eq:interpretvarphiZ}, the mapping $\varphi_Z\left({.\over\varepsilon}\right)$ has a concrete interpretation 
in terms of swaps of points in the Quarck model: For all $i\in\llbracket 1,n-1 \rrbracket$, as the positions of the red points are continuous random variables, we have 
\begin{align}
\varphi_Z\left({B_{(i+1)}-B_{(i)}\over\varepsilon_n}\right)
                                              &= \P\left(Z_i-Z_{i+1}\ge {1\over\varepsilon_n}\left(B_{(i+1)}-B_{(i)}\right) \mid B_{(i+1)}-B_{(i)}\right)\notag\\
                                              &=\P\left(B_{(i)}+\varepsilon_n Z_i \ge B_{(i+1)}+\varepsilon_n Z_{i+1} \mid B_{(i+1)}-B_{(i)}\le \Delta\right)\notag\\
                                              &=\P\left(R_i > R_{i+1} \mid B_{(i+1)}-B_{(i)}\right),\label{eq:representvarphiZ}
\end{align}
namely, the probability of a swap at position $i$ given the $i$-th increment. We deduce the following result.}
\begin{lem}
{Under the ongoing assumptions, if 
$n^2\varepsilon_n\underset{n\to\infty}{\longrightarrow}\infty$ 
we also have 
$\E[V_n]\underset{n\to\infty}{\longrightarrow}\infty.$}
\label{lemma: criticalscaleP1}
\end{lem}
\begin{proof}
Recall \eqref{eq:representvarphiZ}. Applying the law of total probability, we write 
\begin{align*}
 \E[V_n] &=\sum_{i=1}^{n-1}\P\left(R_i>R_{i+1}\right)\\
 &=\sum_{i=1}^{n-1} \P\left(B_{(i)}+\varepsilon_n Z_i \ge B_{(i+1)}+\varepsilon_n Z_{i+1}\right)\\
 &=\sum_{i=1}^{n-1} \P\left(Z_i-Z_{i+1}\ge {1\over\varepsilon_n}\left(B_{(i+1)}-B_{(i)}\right)\right)\\
 &= \sum_{i=1}^{n-1}\sum_{j=1}^{n}\sum_{k=1;k\ne j}^n\P\left(Z_i-Z_{i+1}\ge {1\over\varepsilon_n}\left(B_{k}-B_{j}\right) \mid (i)=j;(i+1)=k\right)\P\left((i)=j;(i+1)=k\right)\\
 &= {1\over n(n-1)}\sum_{i=1}^{n-1}\sum_{j=1}^{n}\sum_{k=1;k\ne j}^n\P\left(Z_i-Z_{i+1}\ge {1\over\varepsilon_n}\left(B_{k}-B_{j}\right) \mid (i)=j;(i+1)=k;B_j < B_k\right)\\
 &= {1\over n(n-1)}\sum_{i=1}^{n-1}\sum_{j=1}^{n}\sum_{k=1;k\ne j}^n
 \frac{\P\left(Z_i-Z_{i+1}\ge {1\over\varepsilon_n}\left(B_{k}-B_{j}\right)\,;\,(i+1)=k\mid (i)=j;B_j < B_k\right)}{\P\left((i+1)=k\mid (i)=j;B_j < B_k\right)}\\
 &= {1\over n(n-1)}\sum_{i=1}^{n-1}\sum_{j=1}^{n}(n-j)\sum_{k=1;k\ne j}^n
 \P\left(Z_i-Z_{i+1}\ge {1\over\varepsilon_n}\left(B_{k}-B_{j}\right)\,;\,(i+1)=k\mid (i)=j;B_j < B_k\right)\\
 &\ge {1\over n(n-1)}\sum_{i=1}^{n-1}\sum_{j=1}^{n}(n-j)\sum_{k=1;k\ne j}^n\P\left(Z_i-Z_{i+1}\ge {1\over\varepsilon_n}\left(B_{k}-B_{j}\right) \mid B_j < B_k\right)\\
 &= {1\over n(n-1)}\sum_{i=1}^{n-1}\sum_{j=1}^{n}(n-j)\sum_{k=1;k\ne j}^n\P\left(Z_j-Z_{k}\ge {1\over\varepsilon_n}\left(B_{k}-B_{j}\right) \mid B_j < B_k\right)\\
 &={(n-1)(n-1)\over n(n-1)}\P\left(Z_1-Z_{2}\ge {1\over\varepsilon_n}\left(B_{2}-B_{1}\right) \mid B_1 < B_2\right)\sum_{j=1}^{n}(n-j)\\
 &={(n-1)^2\over 2}\P\left(Z_1-Z_{2}\ge {1\over\varepsilon_n}\left(B_{2}-B_{1}\right) \mid B_1 < B_2\right),
\end{align*} 
using the independence of the $Z_j-Z_k$'s, for all $j,k$ with the $B_i$'s for all $i$ in the sixth equality, and symmetry in the seventh. But, again by symmetry we have  
\begin{align*}
\P\left(Z_1-Z_{2}\ge {1\over\varepsilon_n}\left(B_{2}-B_{1}\right) \mid B_1 < B_2\right) 
&= 2\P\left(\left\{Z_1-Z_{2}\ge {1\over\varepsilon_n}\left(B_{2}-B_{1}\right)\right\} \cap \{B_1 < B_2\}\right)\\
&=2\int_{[0,1]}\int_{[x,1]}\varphi_Z\bigg(\frac{y-x}{\varepsilon_n}\bigg)f_B(y)f_B(x)\dd\LebR(y)\dd\LebR(x)\\
&=2\int_{[0,1]}f_B(x)\int_{[0,1-x]}\varphi_Z\bigg(\frac{\Delta}{\varepsilon_n}\bigg)f_B(x+\Delta)\dd\LebR(\Delta)\dd\LebR(x)\\
&=2\varepsilon_n\int_{[0,1]}f_B(x)\int_{[0,{1-x\over \varepsilon_n}]}\varphi_Z(\delta)f_B(x+\varepsilon_n\delta)\dd\LebR(\delta)\dd\LebR(x)
\end{align*}
The right factor is positive by assumption. Thus if $n^2\varepsilon_n\rightarrow\infty$ as $n\rightarrow\infty$ so does $\E[V_n]$\,.

\end{proof}
\begin{remark}
    If we have equality in the first equation, we would expect that at the critical scale we have
\[\E[V_n]\xrightarrow{n\rightarrow\infty}\frac\alpha2\E[|Z_1-Z_2|]\norm{f_B}^2_2.\]
\end{remark}
\begin{prop}[\thispropname]
    \label{prop.f+Z->transpositions}
    Let $B_1,\dots,B_n$ be distributed according to a probability density $f_B$ with support contained in the unit interval $[0,1]$ and $\norm{f_B}_p<\infty$ for $p\geq2$. Let the shifts $\oneton Z$ be such that $\E[|Z_1|^k]<\infty$ for some $k>3/2$. 
    Then, for any $\rho>0$, we have 
    \[\sum_{i,k\in\llbracket 1,n\rrbracket:\,i+1<k}\!\P(R_i>R_k)=\OO\left(n^{-\min\left(2(p-1)/p+1;2k-3-\rho\right)}\right).\]
\end{prop}
\begin{proof}
    If we pick some interval $A\subseteq[0,1]$, we have, by Jensen's inequality,
    \[\infty>\int_A f_B^p\dd\LebR\geq\LebR(A)\bigg(\frac1{\LebR(A)}\int_Af_B\dd\LebR\!\bigg)^{\!p}=\LebR(A)^{1-p}\bigg(\int_Af_B\dd\LebR\!\bigg)^p\!.\]
    Hence for a sequence of measurable sets $[0,1]\supset A_1\supset A_2\supset A_3\supset\dots$ with $\lim_{n\rightarrow\infty}\lambda(A_n)=0$, we obtain
    \[\lambda(A_n)^{(1-p)/p}\int_{A_i}f_B\dd\lambda\leq\bigg(\int_{A_n} f_B^p\dd\LebR\bigg)^{1/p}\xrightarrow{n\rightarrow\infty}0.\]
    Thus,
    \[\int_{A_n}f_B\dd\LebR=\oo(\lambda(A_n)^{(p-1)/p}).\] Denote for all $n$, 
    \begin{align*}
        W_n&\coloneq\sum_{i,k\in\llbracket1,n\rrbracket,\,i+1<k}\!\P(R_i>R_k)\\
        &=\E[\sum_{i,k\in\llbracket1,n\rrbracket,\,i+1<k}\bfo[B_{(i)}+\varepsilon_n Z_i>B_{(k)}+\varepsilon_n Z_k]],
    \end{align*}
    and have, as we demanded $i+1<k$ there is a blue point at $x$ between the $i$-th and $k$-th blue points,
    \[W_n=n\int_{[0,1]}f_B(x)\sum_{\overset{i\neq k}{B_i\leq x\leq B_k}}\P(B_i+\varepsilon_n Z_i>B_k+\varepsilon_n Z_k)\dd\LebR(x).\]
    As in the proof of Lemma \ref{lemma.intfxfx+y->ltwonorm_f}, we trivially extend $f_B$ outside of $[0,1]$ by defining $\tilde f_B=f_B\mathbf 1_{[0,1]}$. 
    In the last sum we have $(n-1)(n-2)/2$ choices for $i,k$, thus 
    \[W_n\leq\OO(n^3)\int_{[0,1]}\!\!\varphi_Z\bigg(\frac\Delta{\varepsilon_n}\bigg)\int_{[0,1]}\tilde f_B(x)\int_{[x-\Delta,x]}\tilde f_B(y)\tilde f_B(y+\Delta)\dd\LebR(y)\dd\LebR(x)\dd\LebR(\Delta).\]
    We use \cref{lemma.intfxfx+y->ltwonorm_f} to obtain
    \begin{align*}
        W_n&\leq\OO(n^3)\varepsilon_n\int_{[0,1/\varepsilon_n]}\!\varphi_Z(\Delta)\int_{[0,1]}\tilde f_B(x)\int_{[x-\varepsilon_n\Delta,x]}\tilde f_B(y)\tilde f_B(y+\varepsilon_n\Delta)\dd\LebR(y)\dd\LebR(x)\dd\LebR(\Delta)\\
        &=\OO(n^3)\varepsilon_n\int_{[0,1/\varepsilon_n]}\!\varphi_Z(\Delta)\int_{[0,1-\varepsilon_n\Delta]}\!\tilde f_B(y)\tilde f_B(y+\varepsilon_n\Delta)\!\int_{[y,y+\varepsilon_n\Delta]}\!\tilde f_B(x)\dd\LebR(x)\dd\LebR(y)\dd\LebR(\Delta)\\
        &=\OO(n^3)\varepsilon_n\int_{[0,1/\varepsilon_n]}\!\varphi_Z(\Delta)\int_{[0,1-\varepsilon_n\Delta]}\tilde f_B(y)\tilde f_B(y+\varepsilon_n\Delta)\oo((\varepsilon_n\Delta)^{(p-1)/p})\dd\LebR(y)\dd\LebR(\Delta)\\
        &\leq\OO(n^3)\varepsilon_n\ltwonorm{f_B}^2\int_{[0,1/\varepsilon_n]}\!\varphi_Z(\Delta)\oo((\varepsilon_n\Delta)^{(p-1)/p})\dd\LebR(\Delta).
    \end{align*}
    This gives, for $\nu\geq0$,
    \begin{align*}
        W_n&\leq\OO(n^3)\oo(\varepsilon_n^{3/2})\OO(\varepsilon_n^\nu)\ltwonorm{f_B}^2\int_{[0,1/\varepsilon_n]}\!\varphi_Z(\Delta)\OO(\varepsilon_n^{(p-1)/p-\nu-1/2}\Delta^{(p-1)/p})\dd\LebR(\Delta)\\
         &\leq\OO(n^3)\oo(\varepsilon_n^{3/2})\OO(\varepsilon_n^\nu)\ltwonorm{f_B}^2\cdot\int_{[0,1/\varepsilon_n]}\!\varphi_Z(\Delta)\OO(\Delta^{\nu+1/2})\dd\LebR(\Delta).
    \end{align*}
     Notice that the second inequality only works for $\nu\leq (p-1)/p-1/2$. We remind ourselves that $\varphi_Z(\Delta)=\P(Z_1-Z_2>\Delta)$. Hence, since the $k$-th moment is finite, by the Markov inequality, $\varphi_Z(\Delta)$ is of order $\OO(\Delta^{-k})$. Thus the right factor converges to some finite value if the remaining term is of lower order than $\OO(\Delta^{k-1})$. The right factor converges to some constant as $n\rightarrow\infty$ if $\nu<k-3/2$. The left factor converges to $0$ as $n\rightarrow\infty$ if $\ltwonorm{f_B}<\infty$.
\end{proof}
\section{Critical scale: Proof of Theorem \ref{thm:main}}
\label{sec:proofthm}
In addition to the assumptions of Section \ref{sec:prelimres}, let us now assume that 
the Quarck model $\QM{B}{Z}{\varepsilon_n}{n}$ is at the critical scale. 
To prepare the way for the proof of Theorem \ref{thm:main}, let us first observe the following asymptotics. 
\begin{lem}\label{lemma:convergence.MAD*l2norm}
    {At the critical scale, we have the convergence 
    \[n(n-1)\varepsilon_n\int_{[0,{1/\varepsilon_n}]}\!\!\varphi_Z(\Delta)\int_{[0,{1-\varepsilon_n\Delta}]}\!f_B(x)f_B(x+\varepsilon_n\Delta)\dd\LebR(x)\dd\LebR(\Delta)
    \underset{n\to\infty}{\longrightarrow}\frac\alpha2\E[|Z_1-Z_2|]\ltwonorm{f_B}^2.\]}
\end{lem}
\begin{proof}
    Choose an arbitrarily small $\rho\in(0,\E[|Z_1-Z_2|]\ltwonorm{f_B}^2)$ and use \cref{lemma.intfxfx+y->ltwonorm_f} to take $0<y<1$ such that, for all $0<z\leq y$,
    \[\ltwonorm{f_B}^2-\int_{[0,{1-z}]}\!f_B(x)f_B(x+z)\dd\LebR(x)<\frac\rho{\E[|Z_1-Z_2|]}\cdot\]
    Further choose $n$ large enough ($\varepsilon_n=\alpha n\invv$), using monotone convergence from below and \cref{lemma.MAD}, such that
    \[\frac12\E[|Z_1-Z_2|]-\int_{[0,{y/\varepsilon_n}]}\varphi_Z(\Delta)\dd\LebR(\Delta)<\frac\rho{2\ltwonorm{f_B}^2}\cdot\]
    This gives
    \begin{multline*}
        \int_{[0,{1/\varepsilon_n}]}\!\!\varphi_Z(\Delta)\int_{[0,{1-\varepsilon_n\Delta}]}\!\!f_B(x)f_B(x+\varepsilon_n\Delta)\dd\LebR(x)\dd\LebR(\Delta)\\
        \begin{aligned}
        &\geq\int_{[0,{y/\varepsilon_n}]}\!\!\varphi_Z(\Delta)\int_{[0,{1-\varepsilon_n\Delta}]}\!\!f_B(x)f_B(x+\varepsilon_n\Delta)\dd\LebR(x)\dd\LebR(\Delta)+\int_{[{y/\varepsilon_n},{1/\varepsilon_n}]}\!\!\varphi_Z(\Delta)\cdot0\dd\LebR(\Delta)\\
        &\geq\int_{[0,{y/\varepsilon_n}]}\!\!\varphi_Z(\Delta)\bigg(\!\ltwonorm{f_B}^2-\frac\rho{\E[|Z_1-Z_2|]}\bigg)\dd\LebR(\Delta)\\
        &=\bigg(\!\ltwonorm{f_B}^2-\frac\rho{\E[|Z_1-Z_2|]}\bigg)\int_{[0,{y/\varepsilon_n}]}\varphi_Z(\Delta)\dd\LebR(\Delta)\\
        &\geq\bigg(\!\ltwonorm{f_B}^2-\frac\rho{\E[|Z_1-Z_2|]}\bigg)\bigg(\frac12\E[|Z_1-Z_2|]-\frac\rho{2\ltwonorm{f_B}^2}\bigg)\\
        &>\frac12\E[|Z_1-Z_2|]\ltwonorm{f_B}^2-\rho\cdot
        \end{aligned}
    \end{multline*}
    Conversely, using the Cauchy-Schwarz inequality,
    \begin{align*}
        \int_{[0,{1-\varepsilon_n\Delta}]}\!\!f_B(x)f_B(x+\varepsilon_n\Delta)\dd\LebR(x)&\leq\sqrt{\int_{[0,{1-\varepsilon_n\Delta}]}\!\!f_B^2(x)\dd\LebR(x)\int_{[0,{1-\varepsilon_n\Delta}]}\!\!f_B^2(x+\varepsilon_n\Delta)\dd\LebR(x)}\\&\leq\ltwonorm{f_B}^2,
    \end{align*}
    hence,
    \begin{align*}
        \int_{[0,{1/\varepsilon_n}]}\!\!\varphi_Z(\Delta)\int_{[0,{1-\varepsilon_n\Delta}]}\!\!f_B(x)f_B(x+\varepsilon_n\Delta)\dd\LebR(x)\dd\LebR(\Delta)
        &\leq\ltwonorm{f_B}^2\int_{[0,{1/\varepsilon_n}]}\!\!\varphi_Z(\Delta)\dd\LebR(\Delta)\\
        &\leq\frac12\E[|Z_1-Z_2|]\ltwonorm{f_B}^2.
    \end{align*}
\end{proof}
\noindent We are in a position to prove our main result, Theorem \ref{thm:main}. 
\begin{proof}[Proof of Theorem \ref{thm:main}]
We use the method of factorial moments. In this context, it suffices to show that, for all $k\in\N$,
\[\E[\fallfac{k}{V_n}]\xrightarrow{n\rightarrow\infty}\mathfrak{c}_\alpha^k.\]
    We use the fact that, for a sum of indicator functions, $V_n=\sum_{i=1}^{n-1}\bfo[R_i>R_{i+1}]$, the factorial moment can be written as
\[\E[\fallfac{k}{V_n}]=\sum_{i_1,\dots,i_k\in\llbracket1,n-1\rrbracket\text{ pairwise distinct}}\E[
\bfo[R_{i_1}>R_{i_1+1}]\dots \bfo[R_{i_k}>R_{i_k+1}]],\] together with Proposition \ref{prop.f+Z->transpositions}, to obtain
    \begin{align*}
        \E[\fallfac{k}{V_n}]&=\E[\sum_{i_1\neq\dots\neq i_k\in\llbracket1,n\rrbracket}\prod_{j=1}^k\bfo[R_{i_j}>R_{i_j+1}]]\\
        &=\E[\sum_{\overset{i_1,\dots,i_k,\tilde{i_1},\dots,\tilde{i_k}\in\llbracket1,n\rrbracket}{i_1\neq i_2\neq\dots\neq i_k,\forall j:i_j<\tilde i_j}}\prod_{j=1}^k\bfo[R_{i_j}>R_{\tilde{i_j}}]]+\oo(1)\\
        &=\frac{n!}{(n-2k)!}\E[\prod_{j=1}^k\bfo[\!B_{i_j}+\varepsilon_n Z_{\sigma({i_j})}>B_{\tilde i_j}+\varepsilon_n Z_{\sigma(\tilde i_j)}]\bfo[B_{i_j}<B_{\tilde i_j}]]+\oo(1).
    \end{align*}
    Since the $B_{i_j}$ are independent of each other, we have, by \cref{lemma:convergence.MAD*l2norm}, that
    \begin{align*}
        \E[\fallfac{k}{V_n}]&=\frac{n!}{(n-2k)!}\Bigg(\int_{[0,1]}\varphi\bigg(\frac\Delta{\varepsilon_n}\bigg)\int_{[0,{1-\Delta}]}\!f_B(x)f_B(x+\Delta)\dd\LebR(x)\dd\LebR(\Delta)\Bigg)^k\!+\oo(1)\\
        &=\frac{n!}{(n-2k)!}\bigg(\frac\alpha2n\invv\E[|Z_1-Z_2|]\ltwonorm{f_B}^2\bigg)^k+\oo(1)\\
        &=\bigg(\frac\alpha2\E[|Z_1-Z_2|]\ltwonorm{f_B}^2\bigg)^k+\oo(1),
    \end{align*}
    where for the first line we ignore all cases where some $i_j=\tilde i_{j'},j\neq j'$ since then, assuming $B_{i_j}+\varepsilon_n Z_{\sigma({i_j})}>B_{\tilde i_j}+\varepsilon_n Z_{\sigma(\tilde i_j)},B_{i_{j'}}+\varepsilon_n Z_{\sigma({i_{j'}})}>B_{\tilde i_{j'}}+\varepsilon_n Z_{\sigma(\tilde i_{j'})}$, we also have  $B_{i_{j'}}+\varepsilon_n Z_{\sigma({i_{j'}})}>B_{\tilde i_j}+\varepsilon_n Z_{\sigma(\tilde i_j)}$ whose probability vanishes by \cref{prop.f+Z->transpositions}. We conclude using the method of factorial moments.
    \end{proof}

\begin{cor}
    Let the Quarck model $\QM{B}{Z}{\varepsilon_n}{n}$ be with a density $f_B$ for the blue points with $\norm{f_B}_p<\infty$ for some $p>2$ and let the $k$-th moment for $Z$ exist for some $k>3/2$. We have, as $n\rightarrow\infty$,
    \[\E[V_n]\rightarrow0,\qquad\text{if }\varepsilon_n\in\oo(n\invv).\]
\end{cor}

\section{Application to satellite data}
\subsection{Framework}
Given two consecutive images of the satellite of the same area with some objects of interest marked in the first image, the blue points, and the same points marked in the second image, red points, we assume that the Quarck model describes the situation well. In this situation, the assignment produced by solving the ERAP is hard to verify as the correct assignment is unknown. 

In this section we describe methods to estimate the parameter of the Poissonian limit, including its confidence interval, at the critical scale with and without knowing the correct assignment. This not only yields the possibility of comparing the theoretic results with synthetic distributions and shifts induced by real terrain but also allows us to estimate the number of errors made in practical application without knowledge of the correct assignment.

We observe sets of $n$ blue points $\oneton b$ and $n$ red points $\oneton r$. The distribution of the blue points $f_B$ and the distribution of the shifts are unknown. To provide an easy way to calculate a confidence interval, we may assume that we can gather this data from $k$ independent identical trials.

One way to estimate $\ltwonorm{f_B}^2$ is,
\[\widehat{\ltwonorm{f_B}^2}=\frac{\log n}{n-1}\sum_{1\leq i<j\leq n}\bfo[(n/\log n)|b_i-b_j|\leq1],\]
if $f_B$ is bounded and for some $1/4<\alpha\leq1/2$,
\[\int_\R\int_\R\frac{(f(x-t)-f(x))^2}{|t|^{1+2\alpha}}\dd\LebR(x)\dd\LebR(t)<\infty\]
holds \cite{l2_1}\cite[Theorem 1]{l2_2}. To be able to obtain a reasonably accurate confidence interval, $(\log n)^{2\alpha} n^{-2\alpha+1/2}\ll k^{-1/2}$ shall be true. For further details, we refer to the appendix.
\paragraph{}
We also want to estimate the mean absolute error $\theta=\E[|Z_1-Z_2|]$\,. We assume the assignment obtained by solving the ERAP in one dimension is the correct assignment. We calculate with the observed distances $z_1=r_{(1)}-b_{(1)},\dots,z_n=r_{(n)}-b_{(n)}$. This is justified by the fact that for increasing $n$, the number of transpositions in the assignment has bounded expectation and the deviation induced by a transposition itself is bounded by $1$. Thus, the relative impact approaches $0$ as $n\rightarrow\infty$. We calculate
\[\widehat\MAD=\frac1{n^2}\sum_{i=1}^n\sum_{j=1}^n|z_i-z_j|.\numberthis\label{eq.MAE}\]
This is computationally expensive for large numbers of points as it takes $\OO(n^2)$ time. Hence, we order the $z_i$ into $z_{(1)},\dots,z_{(n)}$ and have, since $n-i$ of the $z_j$ are bigger than $z_{(i)}$, $2((n-i)(-z_{(i)})+(i-1)z_{(i)})$ for each $z_{(i)}$ respectively. Thus,
\[\widehat\MAD=\frac2{n^2}\sum_{i=1}^n(2i-n-1)z_{(i)}.\]
To obtain a confidence interval on $\hat\theta$, we need to compute
\begin{align*}
    \widehat\Var(\widehat\MAD)&=\frac1{n^2-1}\sum_{i=1}^n\sum_{j=1}^n (|z_i-z_j|+\hat\theta)^2\\
    &=\frac1{n^2-1}\sum_{i=1}^n\sum_{j=1}^n ((z_i-z_j)^2-2\hat\theta|z_i-z_j|+\hat\theta^2)\\
    &=\frac1{n^2-1}\sum_{i=1}^n\sum_{j=1}^n(z_i-z_j)^2-2\hat\theta\frac1{n^2-1}\sum_{i=1}^n\sum_{j=1}^n|z_i-z_j|+\frac{n^2}{n^2-1}\hat\theta^2\\
    \text{use Equation \ref{eq.MAE}}\quad&=\frac1{n^2-1}\sum_{i=1}^n\sum_{j=1}^n(z_i^2-2z_iz_j+z_j^2)-\frac{n^2}{n^2-1}\hat\theta^2\\
    &=\frac{2n}{n^2-1}\sum_{i=1}^n z_i^2-\frac2{n^2-1}\left(\sum_{i=1}^n z_i\right)^2-\frac{n^2}{n^2-1}\hat\theta^2,
\end{align*}
which is possible in linear time; the initial sorting is $\OO(n\log n)$.
\paragraph{}
Ultimately, to obtain the confidence intervals, we use in all cases Student's $t$ distribution quantiles, assuming the variance of the estimator is finite. This in particular requires additionally that the shifts have finite variance.
\paragraph{}
If we have the knowledge about the correct assignment, we can collect the error count $V_n$ directly. Hence, we can get the simulated parameter of the Poisson distribution from the average observed error count. Let $X_1,\dots,X_k$ be \iid random variables known to follow a Poisson distribution with an unknown parameter $\Lambda$. Then $X\coloneq\sum_{i=1}^kX_i$ is Poisson distributed with parameter $\Xi\coloneq k\Lambda$ which we observe as $\hat\Xi$ (\ie the total sum of all observed error counts $V_n$). To obtain a confidence interval, we use the quantiles of the $\chi^2$ distribution \cite{garwood}. We notate $\chi^2_\Xi$ for a $\chi^2$ distribution with parameter $\Xi$ and $Q_\Xi$ for its quantile function. We have the following identity:
\[\sum_{i=0}^{m-1}e^{-\Xi}\frac{\Xi^i}{i!}=1-\P(\chi^2_{2m}\leq2\Xi),\quad m\in\Nn.\]
Let $\Xi_L$ be the lower bound of the confidence interval and $\Xi_U$ the upper bound. Then for a $1-\alpha$ confidence interval one has
\begin{align*}
    \frac\alpha2&=\sum_{i=\Xi}^\infty e^{-\Xi_L}\frac{\Xi_L^i}{i!}=\P(\chi^2_{2\Xi}\leq2\Xi_L),\\
    1-\frac\alpha2&=1-\sum^\Xi_{i=0} e^{-\Xi_U}\frac{\Xi_U^i}{i!}=\P(\chi^2_{2(\Xi+1)}\leq2\Xi_U).
\end{align*}
As a result, we have
\begin{align*}
    \hat\Lambda&=\frac1k\Xi\\
    \hat\Lambda_L&=\frac12\frac1kQ_{2\hat\Xi}\bigg(\frac\alpha2\bigg),\\
    \hat\Lambda_U&=\frac12\frac1kQ_{2(\hat\Xi+1)}\bigg(1-\frac\alpha2\bigg).
\end{align*}

\subsection{First numerical results on real-world data}
\begin{figure}[ht!]
    \centering
    \includegraphics[width=0.5\linewidth]{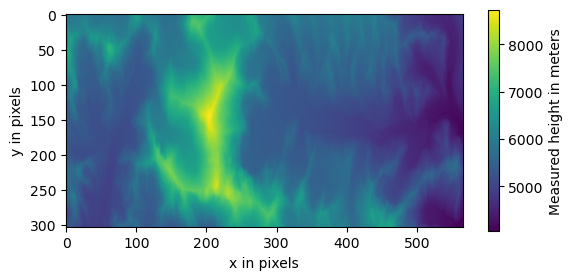}
    \caption[Example image from DEM]{A cropped image from the Copernicus GLO-30 Digital Elevation Model showing Mount Everest and Lhotse.}
    \label{fig:example.image}
\end{figure}

\paragraph{The dataset.}To compare the theoretical results with real-world data, we use the Copernicus GLO-30 Digital Elevation Model \cite{data}. The elevation model is a map of 148.5 million square kilometres of land surface (including inland water bodies) where each pixel represents a $\SI{30}{\metre}\times\SI{30}{\metre}$ square of the surface. We further split the large image into images measuring $566\times303$ pixels; this choice is arbitrary but we choose to be consistent with Quarck's approach \cite{quarck}. An example image showing Mount Everest and Lhotse can be found in \cref{fig:example.image}.
All the data presented in this section was produced using Copernicus WorldDEM-30 \textcopyright DLR e.V. 2010-2014 and \textcopyright Airbus Defence and Space GmbH 2014-2018 provided under COPERNICUS by the European Union and ESA; all rights reserved.

First, we restrict to images that have a mean height of more than \SI{20}{\metre} and a standard deviation of their heights of at least \SI{20}{\metre}. From the remaining data, we extract two sets of data: One set consisting of $3\,221$ images, choosing every 100th available one, effectively scattering them roughly evenly around the globe and one set containing $1\,276$ images containing all images that have a standard deviation of at least \SI{750}{\metre}. The idea of the latter is to include images of very rough terrain where one would expect the shift distribution to have larger $\MAD$. All images with missing data were excluded.
\paragraph{Assumptions.} For the satellite, we choose it to be \SI{200}{\kilo\metre} shifted from the center of the image against the flight direction and it initially moves \SI{200}{\kilo\metre} between taking images. The satellite is assumed to fly \SI{700}{\kilo\metre} above Earth's surface. We neglect distortion effects present in the data and distortion effects that may be produced by the satellite's optics system or the Earth's curvature and in a similar fashion we assume the satellite to fly in a straight line instead of in an orbit. 
To the authors, the distribution of blue points in true observations is unknown. Hence, for all simulations, we assume a uniform distribution and interpolate linearly between pixels. We further assume that no new red points appear and no blue points vanish as they might be optically concealed in real images.

\begin{figure}[thpb!]
    \centering
    \includegraphics[width=0.48\linewidth]{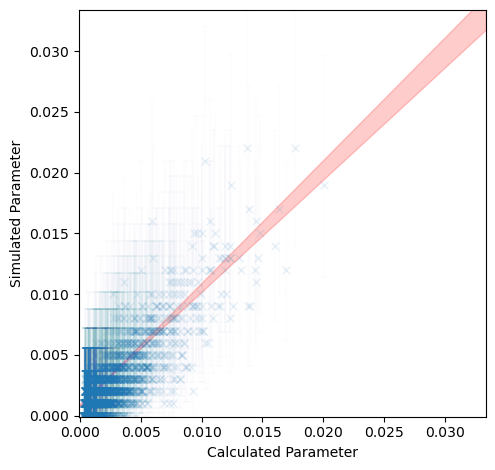}
    \parbox{0.48\linewidth}{\vspace{-6.5cm}\includegraphics[width=\linewidth]{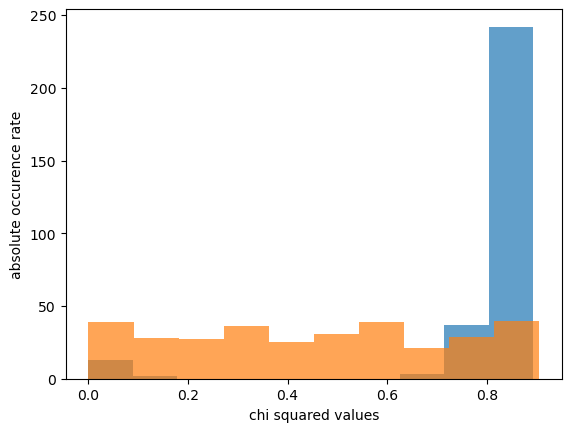}}
    \caption{$1\,000$ simulations with $n=100\,000$ on the normal (first) dataset.
    Left side: Observed calculated parameter plotted against the observed simulated parameter, including the $95\%$ confidence interval obtained from the Monte-Carlo simulations.
    Right side: The results of the $\chi^2$ test for the Poisson distribution for the simulated parameter (blue) and the calculated parameter (orange). Many points are missing since the $\chi^2$ test is not always applicable.}
    \label{fig:spVScp+chi2-normal}
\end{figure}

\begin{figure}[thpb!]
    \centering
    \includegraphics[width=0.48\linewidth]{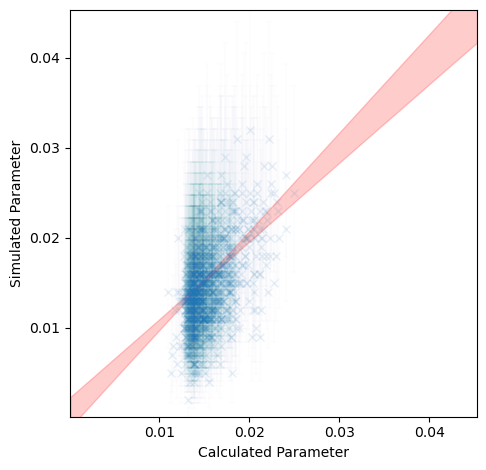}
    \parbox{0.48\linewidth}{\vspace{-6.5cm}\includegraphics[width=\linewidth]{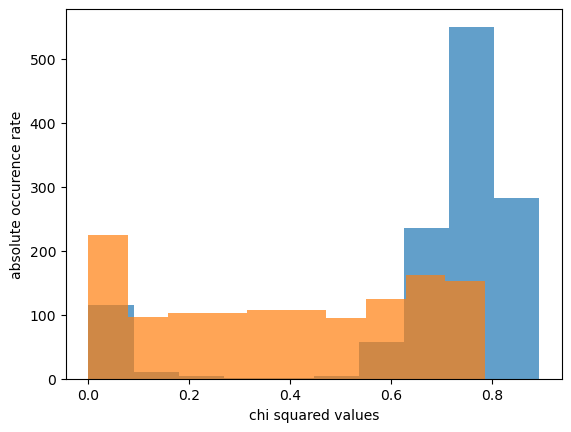}}
    \caption{$1\,000$ simulations with $n=100\,000$ on the dataset with high standard deviation in each image.
    Left side: Observed calculated parameter plotted against the observed simulated parameter, including the $95\%$ confidence interval obtained from the Monte-Carlo simulations.
    Right side: The results of the $\chi^2$ test for the Poisson distribution for the simulated parameter (blue) and the calculated parameter (orange). Some points are missing since the $\chi^2$ test is not always applicable.}
    \label{fig:spVScp+chi2-high}
\end{figure}

\begin{figure}[thpb]
    \centering
    \includegraphics[width=0.48\linewidth]{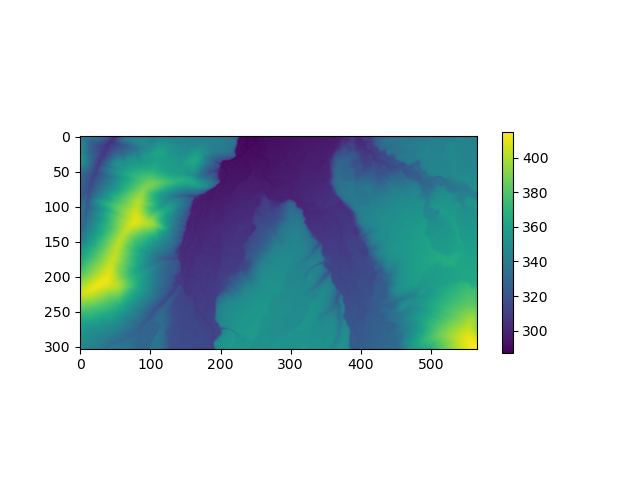}
    \includegraphics[width=0.48\linewidth]{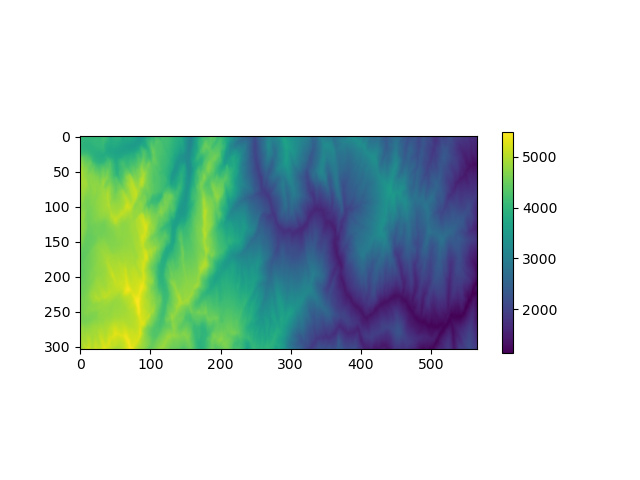}
    \caption{The images used for $500\,000$ simulations. On the left, the image with low standard deviation in its heights. On the right, the image with high standard deviation.}
    \label{fig:examples}
\end{figure}

We are interested in the correspondence between the simulated parameter and the parameter predicted by \cref{thm:main}. We can access with the estimator for $\ltwonorm{f_B}^2$ and $\widehat\MAD$ we developed. To this end, we use a Monte-Carlo-based approach sampling the results obtained with each image in the dataset over many iterations drawing from the distributions of the simulated parameter and the calculated estimate observed. For this, we assume that the calculated parameter is normally distributed and the simulated parameter is gamma distributed.
We use the $\chi^2$ test \cite{cochran} to compare the resulting distribution of the error counts $V_n$ over multiple simulations with the mean of the observed error counts $V_n$ and the result obtained with the estimators. We accumulate the tail of the Poisson distribution into one bin such that at least 5 occurrences are expected. This method fails if the parameter of the Poisson distribution is so small that within the number of simulations only the events $0$ and $1$ are likely to be hit. In that case, the number of bins is 2, which leaves 0 degrees of freedom for the $\chi^2$ as the mean is estimated, not given \cite{cramer}.
\paragraph{}
In \cref{fig:spVScp+chi2-normal} and \cref{fig:spVScp+chi2-high} the results of this approach are presented. Most notable is the fact that the calculated estimates seem to be an inconsistent and rather bad match for the observed distribution.
The $95\%$ confidence intervals for the parameter of the linear regression are shown in \cref{tab:linreg}. The expectation is that the intercept is essentially $0$ and the slope is $1$. The result shows that the intercept is likely to be close to $0$. However, the slope might not be, even though we are unable to reject that hypothesis.
\begin{table}[thpb]
    \centering
    \begin{tabular}{l|cc}
         Dataset&Slope ($95\%$ confidence)&Intercept ($95\%$ confidence)\\\hline\hline
         $\sigma>20$&[0.922,1.001]&[0.001,0.001]\\\hline
         $\sigma>750$&[0.871,1.094]&[0.001,0.002]\\\hline
    \end{tabular}
    \caption{The parameters of the confidence intervals of the linear regressions shown in \cref{fig:spVScp+chi2-normal} and \cref{fig:spVScp+chi2-high}.}
    \label{tab:linreg}
\end{table}
To analyse this phenomenon and the fact that the calculated estimate does not match well, we conduct $500\,000$ simulations instead of $1\,000$ on one image from each dataset (\cref{fig:examples}) and also conduct those simulations on another two images randomly drawn from the height distribution of the previous images.
\cref{tab:highlowredrawn} demonstrates that in both scenarios, with high and low standard deviation in the height distribution, the calculated parameter predicted by \cref{thm:main} is not attained, but, instead, the sample mean may be the parameter for a Poisson distribution that emerges anyway. Removing the local correlation between heights by redrawing the image randomly from its height distribution seems to remove this effect largely and the theoretical prediction matches the simulated results more closely. The remaining difference could be a result of the fact that the satellite by itself visually projects points close to each other likely in the same direction, even though--in the redrawn image--no longer likely by the same distance.
\begin{table}[tbhp]
    \centering
    \renewcommand{\arraystretch}{1.5}
    \begin{tabular}{l|c|c|c|c}
         Image&Simulated&$\chi^2$&Calculated&$\chi^2$\\\hline\hline
         $\sigma>\SI{750}{\metre}$&$0.0188\pm0.0004$&0.91&$0.0201\pm3\cdot10^{-6}$&$10^{-8}$\\\hline
         Redrawn, $\sigma>\SI{750}{\metre}$&$0.0134\pm0.0003$&0.18&$0.0136\pm3\cdot10^{-6}$&0.13\\\hline
         $\sigma>\SI{20}{\metre}$&$4.2\cdot10^{-4}\,\substack{+6\cdot10^{-5}\\
             -5\cdot10^{-5}}$&0.98&$4.7\cdot10^{-4}\pm9\cdot10^{-8}$&0.21\\\hline
         Redrawn, $\sigma>\SI{20}{\metre}$&$2.9\cdot10^{-4}\pm5\cdot10^{-5}$&0.99&$3.2\cdot10^{-4}\pm8\cdot10^{-8}$&0.61\\\hline
    \end{tabular}
    \caption{The simulated and calculated parameter after $500\,000$ simulations of $100\,000$ points for an image with high standard deviation $\sigma$, one with low standard deviation and images randomly redrawn from their respective height distributions. We give the $95\%$ confidence intervals for each parameter. For each parameter, we run the $\chi^2$ test for the parameter against the observed distribution.}
    \label{tab:highlowredrawn}
\end{table}
\begin{remark}
    One might be interested in knowing how the result changes when the shifts move ``likely in the same direction''. For instance, assuming $\Cov(Z_i,Z_j)>0$, one would expect the likelihood of a swap to drop. Consider $\Omega=\{1,2,3\}$ where all $\oinO$ have probability $1/3$. Let $Z,Z'$ be \iid and $Y$ be identically distributed with
    \[\begin{matrix}
        Z(1)=0&Z(2)=36&Z(3)=53\\
        Y(1)=36&Y(2)=0&Y(3)=53
    \end{matrix}\]
    One has
    \begin{align*}
        \E[Z]&=\frac{89}{3},\\
        \Cov(Y,Z)&=\frac13\bigg(-\frac{108}3\frac{19}3\cdot2+\frac{49}3\frac{49}3\bigg)=\frac{349}{27}.
    \end{align*}
    However,
    \begin{align*}
        \E[|Z-Z'|]&=\frac19(36+0+53+53+53+17)=\frac{212}9,\\
        \E[|Z-Y|]&=24,
    \end{align*}
    \ie against the initial expectation $\E[|Z-Y|]>\E[|Z-Z'|]$ even though $Y$ and $Z$ are positively correlated. Hence, this situation remains to be investigated.
\end{remark}

\section{Discussion}
This work investigates a previously unknown, impressively simple model, first found by \cite{quarck}, further. We observed that under general assumptions the global behaviour of the error count is dominated by transpositions of adjacent points. This gave rise to a locally contained probability of an error which, in sum, converges in distribution to a Poisson distribution at the critical scale.

The theoretical part is then accompanied by numerical simulations on real data. However, we were able to observe that the mathematical theory does not quite match the simulated results from the images. A potential solution to improve these results would be to consider a model with some underlying non-random or correlated structure emphasizing the projection of the satellite or the local features of terrain. This could perhaps accommodate the data more precisely.
Conversely, our model was able to demonstrate that it fits well enough to make roughly accurate predictions. Improved performance on images redrawn from their height distribution suggests that the theory models a, in regard to the error count, worse scenario than the situation of a real image. Furthermore, we were able to provide ways to estimate the number of errors without access to the `correct' assignment.

Other possible extensions of this work involve, for example, projecting points for the original two-dimensional images in two orthogonal directions and developing a policy on how to combine the results from both assignment problems. 
More broadly, a sharper general solution for the two-dimensional case is desirable. This is especially relevant for noisy data, if noise induces a shift in the orthogonal direction. 

Moreover, terrain geometry, evolution of the system between two snapshots or of the acquisition angle could influence the number of visible points at each time step. A generalization could involve taking into account a different number $N_R$ of red points compared to the number $N_B$ of blue points.

In the more abstract direction, one could also extend to which sizes of permutations occur for different choices of $\varepsilon_n$ and at which rate. There, one could provide a local, scaled description of the occurring permutations in the form of a permuton limit \cite{HOPPEN2013}.

\paragraph{AI statement.}
All of the proofs in this paper were obtained without relying on AI tools. GPT-6.2 Sol Light was used to proofread the manuscript for typographical errors.

\paragraph{Numerical protocol.}
The simulations were run on the TREX supercomputer, owned by the \CNES (CNES). The simulations used 32 cores with \SI{16}{\giga\byte} RAM each for up to 6 hours on an AMD EPYC Milan 7713, 64C (\SI{2.0}{\giga\hertz} - \SI{225}{\watt}) processor.
\paragraph{Acknowledgements.}
The work of M.D'A.~is partially supported by the LOUCCOUM project of the French Research National Agency (ANR-24-CE40-7809). This work was funded by the Centre National d'Études Spatiales (CNES) under procurement contract No.~5700014405 (D'Achille).
The authors thank S. Ankirchner for useful comments on an earlier version of this manuscript and~N. Gasnier for providing data, and for his guidance and feedback on the application to real-world data.
J.L. has been supported by a postdoctoral grant for CNES. J.L. is grateful to the Institut Élie Cartan de Lorraine (IECL) in Metz for excellent working conditions on
the occasion of two invitations (April and September 2026), during which this work was partly done.
A.S. is grateful to CNES for kind hospitality in the summer of 2026 when part of this work was done.

\bibliographystyle{alpha}

\bibliography{biblio}

\appendix
\section{Estimator for \texorpdfstring{$L^2$}{L2}-norm}\label{app:l2-norm}
To estimate $\norm{f_B}_2^2$, we use the approach of a kernel estimator $T_n(h_n)$ from \cite{l2_1}. We follow the notations of and the results from Gin\'e and Nickl \cite{l2_2}. The estimator is given by
\[T_n(h_n)\coloneq\frac2{h_nn(n-1)}\sum_{1\leq i<j\leq n}K\bigg(\frac{b_i-b_j}{h_n}\bigg),\]
where $K$ is some symmetric bounded kernel function such that $K$ is a probability density with $\int_\R|K(x)||x|\dd\LebR(x)<\infty$. Further one has to choose a sequence $h_n>0$ for which $h_n\rightarrow0$ as $n\rightarrow\infty$.

Since we are also interested in a confidence interval from $k$ trials, the bias of $T_n(h_n)$ has to vanish faster than the standard deviation. We calculate the order of $\Var(T_n(h_n))$. Notice 
\begin{align*}
    \E[\bigg(\frac1n\sum_{i=1}^nY_i\bigg)^2]&=\frac1{n^2}\sum_{i=1}^n\sum_{j=1}^n\E[Y_iY_j]\\
    &=\frac1{n^2}\bigg(n\E[Y_1^2]+\binom n2\E[Y_1Y_2]\bigg)\\
    &=\frac4n\big(\E[f^2(b_1)]-\ltwonorm{f_B}^4\big)+\frac{n-1}{2n}\E[Y_1]^2\\
    &=\frac{4\tau^2}n+0.
\end{align*}
We can obtain using the Cauchy–Bunyakovsky–Schwarz inequality
\begin{align*}
    &\E[\bigg(T_n(h_n)-\E[T_n(h_n)]-\frac1n\sum_{i=1}^nY_i\bigg)^2]\\=&\Var(T_n(h_n))-2\E[\bigg(T_n(h_n)-\E[T_n(h_n)]-\frac1n\sum_{i=1}^nY_i\bigg)\bigg(\frac1n\sum_{i=1}^nY_i\bigg)]-\E[\bigg(\frac1n\sum_{i=1}^nY_i\bigg)^2]\\
    \leq&\Var(T_n(h_n))-2\sqrt{\E[\bigg(T_n(h_n)-\E[T_n(h_n)]-\frac1n\sum_{i=1}^nY_i\bigg)^2]\E[\bigg(\frac1n\sum_{i=1}^nY_i\bigg)^2]}-\frac{4\tau^2}n\\
    =&\Var(T_n(h_n))-2\sqrt{c_2\max\left\{\frac1{n^2h_n},\frac{h_n^{2\alpha}}n\right\}\OO(n\inv)}-\frac{4\tau^2}n.
\end{align*}
Thus
\[\Var(T_n(h_n))\leq\sqrt{\max\left\{\frac1{n^2h_n},\frac{h_n^{2\alpha}}n\right\}\OO(n\inv)}+\frac{4\tau^2}n+c_2\max\left\{\frac1{n^2h_n},\frac{h_n^{2\alpha}}n\right\}.\]
In contrast, the bias vanishes, with $\OO(h^{2\alpha})$. By $h_n=\OO(n^{-2/(4\alpha+1)})$ the highest polynomial order of $h_n$ -- since $f_B$ is unknown -- can only be $n\inv$ as $\alpha\rightarrow1/4$. However, choosing $h_n$ slightly larger, say $(\log n)/n$, allows $\Var(T_n(h_n))=4\tau^2/n+\oo(1/n)$ in comparison to $\Var(T_n(h_n))=\OO(n\inv)$. This way, we improve the variance of the estimate by a constant. Now coming back to $k$ trials, the standard deviation is known to be of order $1/\sqrt k$. Thus, one has to ensure that
\[\OO(h^{2\alpha})=(\log n)^{2\alpha}n^{-2\alpha}\ll\frac{1}{\sqrt{nk}}\]
to obtain a valid confidence interval.
While not fully satisfactory, some stronger assumption on $\alpha$ seems acceptable in practice.
In our scenario, we choose the kernel $K(x)=1/2\,\bfo[|x|\leq1]$ as well as the before-mentioned partially arbitrary $h_n=(\log n)/n$.
\end{document}